\documentclass[12pt,draft]{article}

\usepackage[cp1251]{inputenc}
\usepackage[T2A]{fontenc}
\usepackage[ukrainian,english,russian]{babel}
\usepackage{amsmath,amsfonts,amssymb}
\usepackage{geometry}
\date{06.08.2008}
\large

\begin{document}
\renewcommand{\abstractname}{}
\fontsize{10.8}{12pt}\selectfont 
\vspace{1cm}

\begin{flushleft}

{\textbf {\large Т. И. Кодлюк, В. А. Михайлец, Н. В. Рева }}

 \vspace{0,5cm}

 \noindent\textbf{\Large Непрерывность по параметру решений}

\textbf{\Large одномерных краевых задач}
\end{flushleft}

\emph{ В работе исследуется непрерывность по параметру решений общих
краевых задач для систем линейных обыкновенных дифференциальных уравнений.
Получено обобщение теоремы Кигурадзе (1987) о корректности таких задач.
Найдены также доcтаточные условия равномерной сходимости матриц Грина к матрице Грина
предельной краевой задачи.}

\vspace{0,5cm}\textbf{1. Введение. }Рассмотрим семейство общих линейных неоднородных
 краевых задач для системы $m\in \mathbb{N}$ дифференциальных
 уравнений первого порядка
 \begin{equation*}
y'(t;\varepsilon)=A(t;\varepsilon)y(t;\varepsilon)+f(t;\varepsilon),\quad
t \in(a,b), \eqno(1_\varepsilon)
\end{equation*}
\begin{equation*}
 U_{\varepsilon}y(\cdot;\varepsilon)= c_{\varepsilon}, \quad
 \varepsilon \in [0,\varepsilon_0].
 \eqno(2_\varepsilon)
  \end{equation*}

Здесь квадратные матрицы-функции $A(\cdot;\varepsilon)\in L(
[a,b]; \mathbb{C}^{m\times
  m})=:L^{m\times m}
  $, вектор-функции  \quad $f(\cdot;\varepsilon)\in L( [a,b];
  \mathbb{C}^{m})=:L^m$,\quad  векторы  \quad
  $c_{\varepsilon}\in\mathbb{C}^{m}$, а линейные непрерывные
  операторы
  $$U_{\varepsilon}:C([a,b];\mathbb{C}^{m})
  \rightarrow\mathbb{C}^{m}.$$

Под решением системы дифференциальных уравнений  $(1_\varepsilon)$
понимается  вектор-функция $y(t;\varepsilon)\in W_{1}^{1}
([a,b];\mathbb{C}^m)=:(AC)^m$, для которой равенство $(1_\varepsilon)$
 выполняется на подмножестве интервала полной меры Лебега, которое может зависить от решения.
 Неоднородное „общее“  краевое условие  $(2_\varepsilon)$ охватывает все
 классические виды краевых условий: задачи Коши, двухточечные  и
 многоточечные, интегральные и смешанные краевые задачи (см. [1] и приведенные там
 ссылки).

  Будем предполагать всюду далее, что выполнено

 \textbf{Предположение $\mathcal{E }$.} \emph{Предельная однородная краевая
 задача}

 \begin{equation*}
y'(t;0)=A(t;0)y(t;0), \quad U_{0}y(\cdot;0)= 0,
  \end{equation*}
\emph{имеет только тривиальное решение.}

Это условие равносильно тому, что неоднородная краевая задача
$(1_0),(2_0)$ имеет  решение при произвольных вектор-функции $f(t;0)\in L^m$ и векторе  $c_{0}\in\mathbb{C}^{m}$.

В работе [2], применительно к случаю вещественнозначных функций, установлена следующая

\textbf{Теорема (И. Т. Кигурадзе).} \emph{ Пусть для задачи}
$(1_0),(2_0)$ \emph{выполнено предположение  $\mathcal{E }$ и при $\varepsilon \rightarrow +0$
следующие условия:}

$ 1)\quad \|A(\cdot;\varepsilon)\|_1 = O(1);$

$ 2)\quad  \|f(\cdot;\varepsilon)\|_1 = O(1);$

$ 3)\quad  \| U_\varepsilon\| = O(1);$

$ 4)\quad  \| \int\limits_a ^t A(s;\varepsilon)ds - \int\limits_a ^t A(s;0)ds\|_\infty \rightarrow 0;$

$ 5)\quad  \| \int\limits_a ^t f(s;\varepsilon)ds - \int\limits_a ^t f(s;0)ds\|_\infty \rightarrow 0;$

$ 6)\quad c_\varepsilon \rightarrow c_0, \quad
\varepsilon\rightarrow +0; $

$ 7)\quad U_\varepsilon y \rightarrow U_0 y, \quad \forall y\in
(AC)^m.$

\noindent \emph{Тогда для достаточно малых } $\varepsilon$
\emph{задача} $(1_\varepsilon),(2_\varepsilon)$ \emph{имеет
единственное решение и}
$$ \|y(\cdot;\varepsilon)-y(\cdot;0)\|_\infty \rightarrow 0,\quad
\varepsilon\rightarrow +0 . $$

Здесь и всюду далее $\|\cdot\|_1$
--- норма в пространстве $L_1 = L,$  а
$\|\cdot\|_\infty$ --- sup-норма.

Примеры показывают, что каждое из условий 1)-7) теоремы является существенным и не может быть отброшено. Доказательство И. Т. Кигурадзе не позволяет исследовать случай комплекснозначных функций, который необходим для ряда приложений.

Цель данной работы --- максимально ослабить условия 1), 4) теоремы
Кигурадзе на поведение семейства матриц-функций $A(\cdot;\varepsilon)$
 в окрестности точки $\varepsilon = 0$, а также найти условия,
которые   обеспечивают равномерную сходимость матриц Грина
рассматриваемых задач к матрице Грина предельной краевой задачи  на квадрате $[a,b]\times [a,b]$ .

\textbf{ 2.
Основные результаты.} Введем в рассмотрение класс $\mathcal{M}^{m}$
параметризованных $\varepsilon$ семейств комплекснозначных матриц-функций
$$R(\cdot;\varepsilon): [0,\varepsilon_0]\rightarrow L^{m\times m},$$
для которых матричное решение  $Z(t;\varepsilon)$ задачи Коши
$$ Z'(t;\varepsilon)= R(t;\varepsilon)Z(t;\varepsilon), \quad Z(a;\varepsilon)\equiv I_m$$
удовлетворяет предельному соотношению
$$\lim\limits_{\varepsilon \rightarrow +0} \|Z(t;\varepsilon) - I_m\|_\infty =0,$$
где $I_m$ -- единичная $(m\times m)$-матрица.

\textbf{Теорема 1.} \emph{В формулировке теоремы Кигурадзе можно
заменить условия  } 1), 4) \emph{ одним более общим условием}
$$ R(t;\varepsilon):= A(t;\varepsilon)- A(t;0) \in
\mathcal{M}^m. \eqno(3)$$

\textbf{Замечание.}  \emph{Условие (3) на коэффициенты системы в
теореме 1 уже нельзя ослабить.
 Оно является необходимым, если}
$U_\varepsilon y\equiv y(a).$

 Положим для матрицы-функции $ R(\cdot)\in L^{m\times
m}$ и вектор-функции $ f(\cdot)\in L^m$:
$$ R^\vee (t):= \int\limits_a ^t R(s)ds,\quad f^\vee (t):= \int\limits_a ^t f(s)ds . $$

\noindent Тогда условия 4) и 5) можно переписать соответственно в
виде:

$4')\quad \|R^\vee (\cdot;\varepsilon)\|_\infty\rightarrow 0,\quad
\varepsilon\rightarrow +0;$

$5') \quad\|f^\vee (\cdot;\varepsilon)-f^\vee
(\cdot;0)\|_\infty \rightarrow 0,\quad \varepsilon\rightarrow +0.$

В работах [ 3 --- 7 ]  найдены  необходимые и достаточные условия
того, что матричная функция  $R(t;\varepsilon)\in \mathcal{M}^m$ при выполнении различных \emph{дополнительных}
предположений.

Примеры из работы [6] показывают, что введенный нами класс $\mathcal{M}^{m}$ не является аддитивным.

Из результатов работ [5, 6] следует, что верно

\textbf{Предложение (А. Ю. Левин).} \emph{Если при $\varepsilon\rightarrow +0$ выполнено любое из четырех условий:}

$(\alpha) \quad \|R(\cdot;\varepsilon)\|_1 = O(1)$

$(\beta) \quad\|R^\vee(\cdot;\varepsilon)R(\cdot;\varepsilon)\|_1
\rightarrow 0,$

$(\gamma) \quad\|R(\cdot;\varepsilon)R^\vee(\cdot;\varepsilon)\|_1
\rightarrow 0,$

$(\Delta)
\quad\|R^\vee(\cdot;\varepsilon)R(\cdot;\varepsilon)-R(\cdot;\varepsilon)R^\vee(\cdot;\varepsilon)\|_1
\rightarrow 0,$

\noindent \emph{то условие (3) равносильно условию} $4')$.

В общем случае
условие  $4')$ не является ни необходимым, ни достаточным для
выполнения условия (3). Ниже будет приведен пример, в котором
выполнено соотношение (3), однако не выполняется ни одно из условий
$(\alpha)$, $(\beta)$, $(\gamma)$, $(\Delta)$ и тем более условие 1)
теоремы Кигурадзе.

Доказательство теоремы 1 выделено в п. 3 данной работы.

 Как известно (см., например [2]) для общей краевой задачи
\begin{equation*}
y'(t)=A(t)y(t)+f(t), \quad (Uy)= 0
  \end{equation*}
\noindent существует матрица Грина, т. е. матричная функция
$G(t,s)\in L_\infty([a,b]\times[a,b];\mathbb{C}^{m\times m})$, с
помощью которой решение задачи может быть представлено в
виде:
$$ y(t)=\int\limits_a
^b G(t,s)f(s) ds,\quad  t\in [a,b],\quad f(\cdot)\in L^m .$$ Эта матрица-функция разрывна на диагонали квадрата. Она может иметь и другие разрывы (см. ниже формулу (14)). Для
семейства краевых задач $(1_\varepsilon),(2_\varepsilon)$
матрица-функция $G(t,s)$ зависит от параметра  $\varepsilon$. Поэтому представляет
интерес вопрос о непрерывности по параметру $\varepsilon$
матричной функции  $G(\cdot, \cdot;\varepsilon)$. Ответ на него дает

 \textbf{Теорема 2. }\emph{ Пусть выполнено предположение  $\mathcal{E }$
 и условия:}

 $ 1)\quad A(\cdot;\varepsilon)-A(\cdot;0)\in \mathcal{M}^m;$

 $ 2)\quad \|U_{\varepsilon} - U_{0}\|\rightarrow 0,\quad \varepsilon\rightarrow +0.$

\emph{Тогда для достаточно
малых }$\varepsilon$ \emph{ существуют матрицы Грина задач }
$(1_\varepsilon),(2_\varepsilon)$ \emph{и равномерно на квадрате} $[a,b]\times [a,b]$

$$ \mathop \|G(t,s;\varepsilon)-G(t,s;0)\|_\infty
\rightarrow 0,\quad \varepsilon\rightarrow +0. $$

Доказательство теоремы 2 выделено в п. 4 данной работы. Там же
приведен пример, который показывает, что условие 2) теоремы 2 нельзя
заменить более слабым  условием \emph{сильной } сходимости
$U_\varepsilon$ к $U_0$, которое равносильно паре условий 3), 7)
теоремы Кигурадзе.

Утверждение теоремы 1 анонсировано без доказательства в [11]. В более слабой форме утверждение теоремы 2 использовано в работах [12, 13] для доказательства равномерной резольвентной аппроксимации операторов Штурма-Лиувилля с сингулярными потенциалами. Подобные дифференциальные операторы встречаются в ряде задач современной математической физики.

{\textbf{ 3. Доказательство теоремы 1}. Сформулируем сначала известное (см., например, [9]) утверждение общего характера,
которое будет многократно  использоваться далее.

Пусть  $\mathcal{A}$ -- банахова алгебра с единицей, а $\emph{Inv}$
$\mathcal{ A}$ -- мультипликативная группа  обратимых элементов
алгебры $\mathcal{A}$.

 \textbf{Лемма 1. }

 (1) \emph{Отображение }$X\mapsto X^{-1}$  \emph{является непрерывным по
 норме }$\mathcal{A}$ \emph{на множестве }$\emph{Inv}$ $\mathcal{ A}$;

 (2) \emph{Отображение }$ (X,Y)\mapsto X\cdot Y$\emph{является непрерывным на }$\mathcal{A}
 \times\mathcal{A}$.

В частности, можно положить в лемме 1
$\mathcal{A}=C([a,b];\mathbb{C}^{m\times m})$.
 Это (некоммутативная при $ m\geq 2$) банахова
алгебра с единицей $I_m$ и нормой $$ \|X\|_C := \mathop {\max}
\limits _{a\leqslant t\leqslant b}|X(t)|, \quad |X|:=\sum \limits
_{i,j} |x_{i,j}|.$$

Можно также положить  $\mathcal{A}=\mathcal{B}([a,b];\mathbb{C}^{m\times m})$. Это более
широкая чем $ (C)^{m\times m}$ (несепарабельная)
банахова алгебра ограниченных матриц-функций  с нормой $ \|\cdot\|_\infty.$

Пусть $Y(t;\varepsilon)$ -- единственное решение матричной краевой
задачи
 \begin{equation*}
Y'(t;\varepsilon)=A(t;\varepsilon)Y(t;\varepsilon),\quad
Y(a;\varepsilon)=I_m.
 \end{equation*}

Отправным моментом в нашем доказательстве теоремы 1 является принцип
редукции А. Ю. Левина [5, 6]. В наших обозначениях он имеет
следующий вид:

\textbf{Лемма 2. }\emph{ Предельное соотношение  }

$$ \|Y(\cdot;\varepsilon)-Y(\cdot;0)\|_\infty \rightarrow 0,\quad
\varepsilon\rightarrow +0  $$ \emph{выполняется в том и только в том
случае, если  }
$$ A(\cdot;\varepsilon)-A(\cdot;0) \in
\mathcal{M}^m.$$

Рассмотрим наряду с исходной неоднородной краевой задачей
$(1_\varepsilon),(2_\varepsilon)$ относительно вектор- функции
$y(t;\varepsilon)$ еще три векторные краевые задачи:
\begin{equation*}
z'(t;\varepsilon)=A(t;\varepsilon)z(t;\varepsilon),\quad
 U_{\varepsilon}z(\cdot;\varepsilon)= c_{\varepsilon},
 \eqno(4_\varepsilon)
  \end{equation*}
\begin{equation*}
x'(t;\varepsilon)=A(t;\varepsilon)x(t;\varepsilon)+f(t;\varepsilon),\quad
x(a;\varepsilon)\equiv 0,
 \eqno(5_\varepsilon)
  \end{equation*}
\begin{equation*}
w'(t;\varepsilon)=A(t;\varepsilon)w(t;\varepsilon)+f(t;\varepsilon),\quad
 U_{\varepsilon}w(\cdot;\varepsilon)\equiv 0.
 \eqno(6_\varepsilon)
  \end{equation*}
Как известно, краевая задача $(5_\varepsilon)$ (задача Коши)
всегда имеет решение и оно единственно.

\textbf{Лемма 3. }\emph{ Если выполнено предположение}
$\mathcal{E}$,  \emph{то каждая из задач:}
$(1_\varepsilon)-(2_\varepsilon),(4_\varepsilon)$ и
$(6_\varepsilon)$  \emph{при достаточно малых значениях параметра}
$\varepsilon$ \emph{ имеет ровно одно решение в классе} $(AC)^m.$

\textbf{Доказательство.} Достаточно  показать, что при малых
$\varepsilon$ однородная краевая задача
\begin{equation*}
y'(t;\varepsilon)=A(t;\varepsilon)y(t;\varepsilon),\quad
 U_{\varepsilon}y(\cdot;\varepsilon)= 0
  \end{equation*}
  имеет только тривиальное решение. Каждое из решений однородного
  дифференциального уравнения имеет вид:
$$y(t;\varepsilon)=Y(t;\varepsilon) \widetilde{c}_\varepsilon,
\quad \widetilde{c}_\varepsilon \in \mathbb{C}^m,$$ где
$Y(t;\varepsilon) $ -- матрицант этого уравнения. Откуда в силу
краевого условия имеем:
$$[U_\varepsilon Y(t;\varepsilon)]\widetilde{c}_\varepsilon\equiv 0,$$
где \emph{i}-тый столбец $(m\times m)$ -- матрицы $[U_\varepsilon
Y(t;\varepsilon)]$ совпадает с действием линейного оператора
$U_\varepsilon$ на \emph{i}-тый столбец матрицы
$Y(t;\varepsilon)$.

Квадратная матрица $[U_\varepsilon Y(t;\varepsilon)]$ непрерывно
зависит от $\varepsilon$ в силу леммы 2 и сильной непрерывности  операторной функции $U_\varepsilon$ при
$\varepsilon =0$. Кроме того, в
силу предположения $\mathcal{E}$
$$\det [U_0
Y(t;0)]\neq 0.$$ Поэтому в некоторой окрестности точки $\varepsilon
=0$  функция
$$\det [U_\varepsilon
Y(t;\varepsilon)]\neq 0.$$ Откуда следует, что в этой окрестности
вектор  $\widetilde{c}_\varepsilon\equiv 0$ и лемма доказана.

Из леммы 3 следует,что при малых $\varepsilon>0
$
$$ y(\cdot;\varepsilon)=z(\cdot;\epsilon)+w(\cdot;\epsilon).$$
Поэтому для доказательства теоремы 1 достаточно показать,что при
ее условиях
\begin{equation*}
\|z(\cdot;\varepsilon)-z(\cdot;0)\|_\infty \rightarrow 0,\quad
\varepsilon\rightarrow +0,
 \eqno(7)
  \end{equation*}
\begin{equation*}
\|w(\cdot;\varepsilon)-w(\cdot;0)\|_\infty \rightarrow 0,\quad
\varepsilon\rightarrow +0.
 \eqno(8)
  \end{equation*}

\textbf{Лемма 4. }\emph{ Пусть выполнены условия теоремы 1. Тогда
справедливо предельное соотношение }(7).

\textbf{Доказательство.} Первое из равенств $(4_\varepsilon)$ дает
нам, что
$$z(t;\varepsilon)=Y(t;\varepsilon) \widetilde{c}_\varepsilon.$$
Откуда в силу второго из равенств $(4_\varepsilon)$ получаем, что
$$[U_\varepsilon Y(t;\varepsilon)]\widetilde{c}_\varepsilon=c_\varepsilon.$$
Поэтому, по уже доказанному, при достаточно малых $\varepsilon> 0$
$$\widetilde{c}_\varepsilon=[U_\varepsilon Y(t;\varepsilon)]^{-1}c_\varepsilon.$$
В силу лемм 1, 2
$$|[U_\varepsilon Y(t;\varepsilon)]^{-1}-[U_0 Y(t;0)]^{-1}|\rightarrow 0,\quad
\varepsilon\rightarrow +0.$$ Кроме того, по условию $
c_{\varepsilon}\rightarrow c_{0}.$ Поэтому  $ \widetilde{c}_{\varepsilon}\rightarrow \widetilde{c}_{0}$ при $
\varepsilon\rightarrow +0$. Откуда
следует нужное нам соотношение (7).

\textbf{Лемма 5. }\emph{ Пусть при $ \varepsilon \rightarrow +0$ выполнены условия:}

$1) \quad A(\cdot;\varepsilon)- A(\cdot;0)\in
\mathcal{M}^m;$

$2) \quad \| f(\cdot;\varepsilon)\|_1= O(1);$

$3) \quad \|f^\vee (\cdot;\varepsilon)-f^\vee
(\cdot;0)\|_\infty \rightarrow 0.$

\noindent \emph{Тогда }
\begin{equation*}
\|x(\cdot;\varepsilon)-x(\cdot;0)\|_\infty \rightarrow 0,\quad
\varepsilon\rightarrow +0.
 \eqno(9)
  \end{equation*}

\textbf{Доказательство.} Из условия 1) в силу принципа редукции
вытекает, что
\begin{equation*} \|Y(t;\varepsilon)-Y(t;0)\|_\infty
\rightarrow 0,\quad \varepsilon\rightarrow +0.
  \end{equation*}
  Откуда в силу леммы 1 следует, что
  \begin{equation*} \|Y^{-1}(t;\varepsilon)-Y^{-1}(t;0)\|_\infty
\rightarrow 0,\quad \varepsilon\rightarrow +0.
  \end{equation*}

Как известно, решение $ x(t;\varepsilon)$ задачи $(5_\varepsilon)$
может быть представлено в виде
$$ x(t;\varepsilon)= Y(t;\varepsilon)\int \limits _a ^t
Y^{-1}(s;\varepsilon)f(s;\varepsilon)ds.$$ Поэтому ввиду леммы 1
достаточно доказать, что
$$\| \int \limits _a ^t
Y^{-1}(s;\varepsilon)f(s;\varepsilon)ds
- \int \limits _a ^t
Y^{-1}(s;0)f(s;0)ds\|_\infty \rightarrow 0.$$
 Из оценки $$\|\int \limits _a
^t [Y^{-1}(s;\varepsilon)-Y^{-1}(s;0)]f(s;\varepsilon)ds\|_\infty
\leqslant \int \limits _a ^t
|Y^{-1}(s;\varepsilon)-Y^{-1}(s;0)|\cdot|f(s;\varepsilon)|ds\leqslant$$

$$\leqslant\|Y^{-1}(s;\varepsilon)-Y^{-1}(s;0)\|_\infty \cdot \sup\limits
_{\varepsilon}\|f(s;\varepsilon)\|_1 \leqslant c \cdot
\|Y^{-1}(s;\varepsilon)-Y^{-1}(s;0)\|_\infty \rightarrow 0$$ вытекает,
что достаточно доказать, что

$$\|\int \limits _a ^t
Y^{-1}(s;0)[f(s;\varepsilon)-f(s;0)]ds\|_\infty \rightarrow 0. $$
Интегрируя интеграл по частям имеем:
$$\|\int \limits _a ^t
Y^{-1}(s;0)[f(s;\varepsilon)-f(s;0)]ds\|_\infty \leqslant $$
$$\leqslant \|\int \limits _a ^t
(Y^{-1})'(s;0)[f^{\vee}(s;\varepsilon)-f^{\vee}(s;0)]ds\|_\infty +
2\|Y^{-1}(s;0)\|_\infty \cdot
\|f^{\vee}(s;\varepsilon)-f^{\vee}(s;0)\|_\infty \leqslant$$ $$
\leqslant\|f^{\vee}(\cdot;\varepsilon)-f^{\vee}(\cdot;0)\|_\infty \cdot
(2\|Y^{-1}(s;0)\|_\infty +\|Y^{-1}(s;0)\|^2_\infty \cdot
\|Y'(s;0)\|_1)\rightarrow 0, \quad \varepsilon\rightarrow +0,
$$
т. к.  $$(Y^{-1})'(\cdot;\varepsilon)=
-Y^{-1}(\cdot;\varepsilon)Y'(\cdot;\varepsilon)Y^{-1}(\cdot;\varepsilon).$$
Лемма доказана.
 \vspace{0,2 cm}

 \textbf{Лемма 6. } \emph{ При условиях теоремы 1
справедливо предельное соотношение }(8).

 \textbf{Доказательство.} Положим
 $$v(t;\varepsilon):=x(t;\varepsilon)-w(t;\varepsilon).$$
Тогда вектор-функция $v(t;\varepsilon)$ является решением
краевой задачи
 $$v'(t;\varepsilon)=A(t;\varepsilon)v(t;\varepsilon),
 \quad U_\varepsilon v(t;\varepsilon)=U_\varepsilon x(t;\varepsilon)
 =:\widetilde{c}_\varepsilon.$$
Но
$$\|U_\varepsilon x(t;\varepsilon)- U_0 x(t;0)\|_C \leqslant
\|U_\varepsilon\|\|x(t;\varepsilon)-x(t;0)\|_C + \|(U_\varepsilon-
U_0) x(t;0)\|_C \rightarrow 0,$$ т. е. $\quad
\widetilde{c}_{\varepsilon}\rightarrow \widetilde{c}_{0}$ при $
\varepsilon\rightarrow +0$. Поэтому

$$v(t;\varepsilon)=Y(t;\varepsilon)\overline{c}_\varepsilon,
\quad\overline{c}_\varepsilon\in\mathbb{C}^m, $$ где
$[U_\varepsilon Y(t;\varepsilon)]\overline{c}_\varepsilon=
\widetilde{c}_{\varepsilon}$ и при достаточно малых $\varepsilon >
0$
$$\overline{c}_\varepsilon=[U_\varepsilon
Y(t;\varepsilon)]^{-1} \widetilde{c}_{\varepsilon}\rightarrow [U_0
Y(t;0)]^{-1} \widetilde{c}_{0}=\overline{c}_0.
$$ Откуда следует, что
\begin{equation*} \|v(t;\varepsilon)-v(t;0)\|_\infty
\rightarrow 0,\quad \varepsilon\rightarrow +0.
 \eqno(10)
  \end{equation*}
Из равенства  $w(t;\varepsilon)=x(t;\varepsilon)-
v(t;\varepsilon)$ и уже доказанных соотношений (9)и (10) следует
асимптотическое соотношение (8).

Лемма 6, а вместе с ней и теорема 1 доказаны.

Приведем пример, в котором выполнено соотношение (3), однако не
выполняется ни одно из условий $(\alpha)$, $(\beta)$, $(\gamma)$,
$(\Delta)$. Он, в частности,  показывает, что теорема 1 сильнее
теоремы Кигурадзе.

 \textbf{Пример 1. } Пусть $m=2$, $\quad (a,b)=(0,1)$,
$\quad A(t;\varepsilon)=A(t)+R(t;\varepsilon)$, где

$$R(t;\varepsilon)=\left(\begin{array}{cc}0&
 {\frac{1}{\sqrt \varepsilon}}\cos{(\frac {t}{\varepsilon})} \\
{\frac{1}{\sqrt \varepsilon}}\sin{(\frac {2t}{\varepsilon})}
&0\end{array}\right).$$

 \noindent Нетрудно проверить, что $\|R^\vee(\cdot;\varepsilon)\|_\infty\rightarrow
 0$ и
$$R(t;\varepsilon)R^\vee(t;\varepsilon)= diag \\\ \{ {\frac {1}{2}}
\sin{(\frac {2t}{\varepsilon})}\cdot\sin{(\frac {t}{\varepsilon})},\quad
\sin{(\frac {2t}{\varepsilon})}\cdot\sin{(\frac {t}{\varepsilon})} \},$$
$$R^\vee(t;\varepsilon)R(t;\varepsilon)= diag \\\ \{\sin{(\frac
{t}{\varepsilon})} \cdot\sin{(\frac {2t}{\varepsilon})},\quad {\frac
{1}{2}}\sin{(\frac {t}{\varepsilon})} \cdot\sin{(\frac
{2t}{\varepsilon})}\},$$
$$R^\vee(t;\varepsilon)R(t;\varepsilon)-R(t;\varepsilon)R^\vee(t;\varepsilon)=
diag \\\ \{ -{\frac {1}{2}} \sin{(\frac
{2t}{\varepsilon})}\cdot\sin{(\frac {t}{\varepsilon})}, \quad {\frac
{1}{2}}\sin{(\frac {t}{\varepsilon})}\cdot\sin{(\frac
{2t}{\varepsilon})}\}.
$$
\noindent Однако:$$\|R(t;\varepsilon)\|_1 \geq
{\frac{1}{\sqrt{\varepsilon}}}\int\limits_0^1
 |\cos{(\frac{t}{\varepsilon})}|dt=\sqrt{\varepsilon}\int\limits_0^{\frac{1}{\varepsilon}}
 |\cos{(t)}|dt={\frac{1}{\sqrt{\varepsilon}}}M\{|\cos{(t)}|\}\rightarrow +\infty,$$
 $$ \int\limits_0^1 |\sin{(\frac {t}{\varepsilon})}\cdot\sin{(\frac
{2t}{\varepsilon})}|dt=\varepsilon
\int\limits_0^{\frac{1}{\varepsilon}}|\sin
{(t)}|\cdot|\sin{(2t)}|dt\rightarrow M\{|\sin{(t)}\cdot\sin{(2t)}|\}>0,$$
[10]. Поэтому ни одно из четырех приведенных выше условий
здесь не выполнено. Однако, пользуясь теоремой 6 работы [7] при
$i=1,$ нетрудно убедиться, что  $ R(\cdot;\varepsilon)\in
\mathcal{M}^2.$

\vspace{0,5cm} \textbf{4. Доказательство теоремы 2.} Рассмотрим полуоднородную векторную  краевую задачу
\begin{equation*}
y'(t)=A(t)y(t)+f(t), \quad Uy= 0, \eqno (11)
  \end{equation*} где $A(\cdot)\in L^{m\times m},$ $f(\cdot)\in L^m$,
   $U\in \mathcal{L}((C)^m ;\mathbb{C}^m)$. Тогда справедливо однозначное представление
$$Uy = \int \limits _a ^b [dH(t)]y(t),\quad y(\cdot)\in (C)^m,$$
где  $H(\cdot)\in NBV([a,b];\mathbb{C}^{m\times
m})$. Это банахово пространство комплекснозначных
$(m\times m)$--матриц-функций с ограниченным изменением на
отрезке $[a,b]$, которые равны $0$ в точке  $a$ и непрерывны слева
на полуинтервале $(a,b]$. Поэтому для матрицанта $Y(t)$ системы (11)
на интервале $[a,b]$ определена заданная интегралом Стильтьеса матрица
- функция:

\begin{equation*}
H_{Y}(t)= \int \limits _a ^t [dH(s)]Y(s).\eqno (12)
\end{equation*}
Она разрывна в точках разрыва матрицы-функции $H(\cdot).$ При этом, если однородная краевая задача (11) имеет только
тривиальное решение, то
\begin{equation*}
\det H_{Y}(b)\neq 0\eqno (13)
\end{equation*} и существует матрица $H^{-1}_Y(b).$

Как и в вещественном случае (см., например, [8]) верна

\textbf{Лемма 7. } \emph{Если выполнено неравенство }(13),
\emph{то матрица Грина задачи} (11) \emph{ существует и представима
в виде:}

$$G(t,s)=\left\{
\begin{array}{lc}
    Y(t)Y^{-1}(s) - Y(t)H^{-1}_Y(b)H_{Y}(s)Y^{-1}(s),  & a \leqslant s\leqslant t\leqslant b; \\
    - Y(t)H^{-1}_{Y}(b)H_{Y}(s)Y^{-1}(s) ,  & a \leqslant t<s\leqslant b. \\
\end{array}
\right.   \eqno (14) $$

Формулу (14) удобно записать в виде $G(t,s) = G_1(t,s) + G_2(t,s)$, где

$$ G_1(t,s) = - Y(t)H^{-1}_Y(b)H_{Y}(s)Y^{-1}(s)$$,
$$ G_2(t,s)=\left\{
\begin{array}{lc}
    Y(t)Y^{-1}(s),  & a \leqslant s\leqslant t\leqslant b; \\
    0 ,  & a \leqslant t<s\leqslant b. \\
\end{array}
\right.$$ Понятно, что достаточно показать, что

$$ \|G_i(t,s;\varepsilon) - G_i(t,s;0)\|_\infty \rightarrow 0, \quad \varepsilon \rightarrow +0,\quad i=1,2.$$ 
Из леммы 1 вытекает, что если
$$\|T_\varepsilon(t) - T_0(t)\|_\infty \rightarrow 0, \quad \|S_\varepsilon(s) - S_0(s)\|_\infty \rightarrow 0, \quad C_\varepsilon \rightarrow C_0,$$
то на квадрате $(t, s) \in [a, b] \times [a, b]$
$$ \|T_\varepsilon(t)C_\varepsilon S_\varepsilon(s) - T_0(t)C_0 S_0(s)\|_\infty \rightarrow 0, \quad \varepsilon \rightarrow +0.$$

В силу этого для доказательства теоремы 2 достаточно
показать, что при выполнении ее условий для достаточно малых
$\varepsilon$
$$ \det H_{Y}(b;\varepsilon)\neq 0,$$
$$\|Y(\cdot;\varepsilon)-Y(\cdot;0)\|_\infty\rightarrow 0,
 \quad \varepsilon\rightarrow +0,$$
$$\|H_{Y}(\cdot;\varepsilon)-H_{Y}(\cdot;0)\|_\infty\rightarrow 0,
 \quad \varepsilon\rightarrow +0.$$ Первое из предельных
 соотношений уже установлено нами.

Переходя ко второму, имеем:

$$\|H_{Y}(t;\varepsilon)-H_{Y}(t;0)\|_\infty=
\|\int \limits _a ^t [dH(s;\varepsilon)]Y(s;\varepsilon)-\int
\limits _a ^t [dH(s;0)]Y(s;0)\|_\infty \leqslant$$ $$
 \leqslant
  \|\int \limits _a ^t [d(H(s;\varepsilon)-H(s;0))]Y(s;\varepsilon)\|_\infty+
  \|\int \limits _a ^t [dH(s;0)]\cdot[Y(s;\varepsilon)-Y(s;0)]\|_\infty\leqslant
 $$ $$ \leqslant \operatorname{Var}_a ^b[H(s;\varepsilon)-H(s;0)]\cdot \|Y(\cdot;\varepsilon)\|_\infty+
 \operatorname{Var}_a ^b[H(s;0)] \cdot\|Y(\cdot;\varepsilon)-Y(\cdot;0)\|_\infty \rightarrow 0,$$
 т.к. в силу условия
$\|U_\varepsilon-U_0\|\rightarrow 0$ вариация матрицы-функции

$$ \operatorname{Var}_a ^b[H(s;\varepsilon)-H(s;0)]\rightarrow 0.$$
В частности, отсюда следует, что $\operatorname{det} H_Y(b;\varepsilon) \neq 0$ для достаточно малых $\varepsilon$.
Теорема доказана.

Приведем обещанный в п.2 пример.

\textbf{Пример 2. } Пусть $m=1,$ $(a,b)=(0,1),$
$A(t,\varepsilon)=A(t,0)=0,$ а линейные непрерывные операторы
$U_\varepsilon: C([0,1];\mathbb{C})\rightarrow \mathbb{C}$ заданы
равенством
$$(U_\varepsilon y):=y(\varepsilon),\quad \varepsilon\in [0,1].$$
Определенные таким образом операторы $U_\varepsilon$ сильно сходятся
к оператору  $U_0$ на пространстве $(C)^m$:
$$(U_\varepsilon y)=y(\varepsilon)\rightarrow y(0),
\quad \varepsilon\rightarrow + 0,\quad  y(\cdot)\in
C([0,1];\mathbb{C}).$$ Однако
$$\|U_\varepsilon-U_0\|=2,\quad \varepsilon\neq 0.$$ В данном
случае функция Грина задачи
$$y'(t;0)=f(t), \quad y(t;0)|_{t=0}=0$$ имеет вид:
$$G(t,s;0)=\left\{%
\begin{array}{ll}
    1, & 0\leqslant s<t\leqslant 1; \\
    0, & 0\leqslant t\leqslant s\leqslant 1;\\
\end{array}%
\right.    $$ а для задачи
$$y'(t;\varepsilon)=f(t),\quad  y(t;\varepsilon)|_{t=\varepsilon}=0,$$
соответственно будет
$$G(t,s;\varepsilon)= G(t,s,0) - \textbf{1}_{[0,\varepsilon]\times [0,\varepsilon]}(t,s), $$ где $\mathbf{1}_F$  -- характеристическая функция множества $F$.

Откуда следует, что
$$ \quad \|G(t,s;\varepsilon)-G(t,s;0)\|_\infty
=1,\quad \varepsilon\neq 0. $$

Для операторов, отвечающим многоточечным краевым задачам с
$$ U_\varepsilon y := B_1(\varepsilon)y(t_1) + B_2(\varepsilon)y(t_2) + \ldots + B_n(\varepsilon)y(t_n),$$ где $n \geq 2$, точки $\{t_1, t_2, \ldots , t_n\} \in [a,b]$ и не зависят от $\varepsilon$, матрицы $B_k(\varepsilon) \in \mathbb{C}^{m\times m}$, условия
$$ \|U_\varepsilon - U_0\| \rightarrow 0, \quad и \quad U_\varepsilon y \rightarrow U_0 y, \quad y \in (AC)^m$$ равносильны между собой. Каждое из них эквивалентно тому, что
$$B_k(\varepsilon) \rightarrow B_k(0), \quad \varepsilon \rightarrow +0, \quad k = 1, 2,\ldots , n.$$

\vspace{0,5 cm}
\begin{enumerate}

\item \textit{Камке Э.} Справочник по обыкновенным
дифференциальным уравнениям. -- Москва: Наука, 1965. -- 703 с.
\item  \textit{Кигурадзе И. Т.} Краевые задачи для систем
обыкновенных дифференциальных уравнений // Совр. проблемы
математики. Новейшие достижения, т. 30. -- Москва: ВИНИТИ, 1987.
-- С. 3-103.

\item\textit{ Reid W. T.}  Some limit theorems for ordinary
differential systems // J. Diff. Equat. -- 1967. -- 3, № 3. -- Р.
423-439.

\item\textit{ Opial Z.}  Continuous parameter dependence in linear
systems of differential equations  // J. Diff. Equat. -- 1967. --
3. -- Р. 571-579.
 \item \textit{ Левин А. Ю.} Предельный переход для несингулярных
 систем $\dot{X}=A_n (t)X$.// Докл. АН СССР. -- 1967. -- 176, № 4.
 -- С. 774-777.

\item \textit{ Левин А. Ю.}  Вопросы теории обыкновенного
линейного дифференциального  уравнения. I. // Вестник Ярославского
университета. -- 1973. -- Вып. 5. -- С. 105-132.

 \item \textit{Нгуен Тхе Хоан}. О зависимости от параметра
решений линейной системы дифференциаль-

ных уравнений // Дифф.
уравнения. -- 1993. -- 29, №6. -- С. 970-975. \item
\textit{Кигурадзе И. Т.}  Некоторые сингулярные краевые задачи для
обыкновенных  дифферен-

циальных уравнений. -- Тбилиси: Из- во
Тбилисского ун-та. -- 1975. -- 352 с.

 \item  \textit{Данфорд Н., Шварц Дж. Т.}
Линейные операторы. Спектральная теория. -- Москва: Мир.-
 -- 1966. -- 1064 с.
 \item \textit{Демидович Б. П.}
Лекции по математической теории
 устойчивости. -- Москва: Наука, 1967. -- 472 с.

 \item \textit{Михайлец В. А.}
 Обобщения теоремы Кигурадзе о корректности линейных краевых задач / В. А. Михайлец, Н. В. Рева //
Доповіді НАН України. -- 2008. -- №~9. -- С.~23--27.

 \item \textit{Горюнов А. С.}
 Резольвентная сходимость операторов Штурма-Лиувилля с сингулярными потенциалами
/ А. С. Горюнов, В. А. Михайлец  //  Математические заметки. -- 2010. -- Т.87,  № 2. -- С. 311--315.

 \item \textit{Goriunov A. S., Mikhailets V. A.}
 Regularization of singular Sturm-Liouville equations / A.~S. Goriunov,  V.~A. Mikhailets  //
Methods of Functional Analysis and Topology. --- 2010. --- V.~16, №~2. --- P. 120---130.
\end{enumerate}

 \newpage

 \end{document}